\documentclass[11pt,bezier]{article}
\usepackage{amsmath}
\usepackage{amsfonts,amsthm,amssymb}
\usepackage{amsfonts}
\usepackage{graphics}
\usepackage{amssymb}
\newtheorem{lem}{Lemma}[section]
\newtheorem{thm}[lem]{Theorem}
\newtheorem{cor}[lem]{Corollary}

\theoremstyle{definition}

\begin{document}
\title{The super restricted edge-connectedness of direct product graphs\footnote{The research is supported by National Natural Science Foundation of China (12261086, 11861066).}}
\author{Jiaqiong Yin, Yingzhi Tian\footnote{Corresponding author. E-mail: tianyzhxj@163.com.} \\
{\small College of Mathematics and System Sciences, Xinjiang
University, Urumqi, Xinjiang, 830046, PR China}}

\date{}
\maketitle

\noindent{\bf Abstract } Let $G$ be a graph with vertex set $V(G)$ and edge set $E(G)$. An edge subset $F\subseteq E(G)$ is called a restricted edge-cut if $G-F$ is disconnected and has no isolated vertices. The restricted edge-connectivity $\lambda'(G)$ of $G$ is the cardinality of a minimum restricted edge-cut of $G$ if it has any; otherwise $\lambda'(G)=+\infty$.  If $G$ is not a star and its order is at least four, then $\lambda'(G)\leq \xi(G)$, where $\xi(G)=$ min$\{d_G(u) + d_G(v)-2:\ uv \in E(G)\}$. The graph $G$ is said to be maximally restricted edge-connected if $\lambda'(G)= \xi(G)$; the graph $G$ is said to be super restricted edge-connected if every minimum restricted edge-cut isolates an edge from $G$. The direct product of graphs $G_1$ and $G_2$, denoted by $G_1\times G_2$, is the graph with vertex set $V(G_1\times G_2) = V(G_1)\times V(G_2)$, where two vertices $(u_{1} ,v_{1} )$ and $(u_{2} ,v_{2} )$ are adjacent in $G_1\times G_2$ if and only if $u_{1}u_{2} \in E(G_1)$ and $v_{1}v_{2} \in E(G_2)$. In this paper, we give a sufficient condition for $G\times K_{n}$ to be super restricted edge-connected, where $K_{n}$ is the complete graph on $n$ vertices.


\section{Introduction}

It is well known that the underlying topology of an interconnection network can be modeled by a graph $G=(V(G), E(G))$, where $V(G)$ represents the nodes of the network and $E(G)$ represents the links between nodes. Under the assumption that the nodes are reliable and the links may be destroyed. The edge-connectivity of the graph $G$ is used to measure the reliability of the corresponding network.

For a vertex $u\in V(G)$, the $neighborhood$ $N_G(u)$ of $u$ is the set of vertices adjacent to $u$  in $G$, and the degree $d_{G}(u)$ of $u$ is $|N_G(v)|$. For an edge $e = uv \in E(G)$, the $edge$-$degree$ $\xi_{G}(e)$ of $e$ is $d_{G} (u)+d_{G}(v)-2$. The $minimum$ $degree$ $\delta(G)$  and the $minimum$ $edge$-$degree$ $\xi(G)$ of $G$ are defined by  min$\{d_{G}(u): u \in V (G)\}$ and min$\{d_G(u) + d_G(v)-2: uv \in E(G)\}$, respectively. A graph is nontrivial if it has at least two vertices.  

An edge subset $F\subseteq E(G)$ is called an $edge$-$cut$ if $G-F$ is not connected. The $edge$-$connectivity$ of $G$, denoted by $\lambda(G)$, is the cardinality of a  minimum edge-cut of $G$. It is well-known that $\lambda(G)\leq \delta(G)$. So $G$ is said to be $maximally$ $edge$-$connected$, or simply $\lambda$-$optimal$, if $\lambda(G)=\delta(G)$. In addition, $G$ is said to be  $super$ $edge$-$connected$, or simply $super$-$\lambda$, if every $minimum$ $edge$-$cut$ isolates a vertex from $G$.

One deficiency of the edge-connectivity is that it only considers when the remaining graph is not connected, but not the properties of the components. For compensating this deficiency, Esfahanian and Hakimi $\cite{Esfahanian&Hakimi}$ introduced the notion of restricted edge-connectivity, which is a kind of conditional edge-connectivity introduced by Harary \cite{Harary}. An edge-cut $F$ is called a $restricted$ $edge$-$cut$ if each component of  $G-F$ has at least two vertices. The $restricted$ $edge$-$connectivity$ of $G$, denoted by $\lambda'(G)$, is the cardinality of a  minimum restricted edge-cut of $G$ if it has any; otherwise $\lambda'(G)=+\infty$. If $G$ is not a star and its order is at least four, then $\lambda^{'}(G)\leq\xi(G)$. So the graph $G$ is called $maximally$ $restricted$ $edge$-$connected$, or simply $\lambda^{'}$-$optimal$, if $\lambda^{'}(G)=\xi(G)$; and the graph $G$ is called $super$ $restricted$ $edge$-$connected$ or simply $super$ -$\lambda^{'}$, if every minimum $restricted$ $edge$-$cut$ isolates an edge from $G$. Clearly, super edge-connected and super restricted edge-connected graphs obtain the maximum edge-connectivity and restricted edge-connectivity among the graphs with given minimum degree and minimum edge-degree, and only have trivial minimum edge-cuts and minimum restricted edge-cuts, respectively.

Let $G_1$ and $G_2$ be two  graphs. The $direct$ $product$ (also named kronecker product, tensor product and cross product et al.) $G_1\times G_2$ has vertex set $V(G_1\times G_2) = V(G_1)\times V(G_2)$ and edge set $E(G_1\times G_2)=\{(u_{1},v_{1})(u_{2},v_{2}): u_{1}u_{2}\in E(G_1)$ and $v_{1}v_{2}\in E(G_2)\}$. Weichsel \cite{Weichsel} proved that the direct product of two nontrivial graphs is connected if and only if both graphs are connected and  at least one of them is not bipartite.

Bre\v{s}ar and  \v{S}pacapan \cite{Bresar} obtained some bounds on the edge-connectivity of the direct product of graphs. Cao,  Brglez,  \v{S}pacapan and  Vumar \cite{Cao} determined the edge-connectivity of direct product of a nontrivial graph and a complete graph. In \cite{Spacapan}, \v{S}pacapan not only obtained the edge-connectivity of direct product of two  graphs, but also characterized the structure of a minimum edge-cut in the direct product of two graphs.

The complete graph on $n$ vertices is denoted by $K_n$. The total graph, denoted by $T_{n}$, is the graph obtained from $K_{n}$ by adding a loop to each vertex of $K_{n}$.
In \cite{Ma},  Ma, Wang and Zhang studied the restricted edge-connectivity of directed product of a nontrivial graph with a complete graph or a total graph.

\begin{thm} (\cite{Ma}) For any nontrivial connected graph $G$ and any integer $n\geq3$,
\[
\lambda^{'}(G\times K_{n})=min\{(n-1)\xi(G)+2(n-2),\ n(n-1)\lambda^{'}(G)\}.
\]
\end{thm}

\begin{thm} (\cite{Ma}) For any nontrivial connected graph $G$ and any integer $n\geq3$,
\[
\lambda^{'}(G\times T_{n})=min\{n\xi(G)+2(n-1),\ n^{2}\lambda^{'}(G)\}.
\]
\end{thm}

In this paper, we study  the super restricted edge-connectedness of the direct product of a nontrivial connected graph with a complete graph or a total graph. In the next section, some notation and preliminary results will be given. The main results will be presented in Section 3.

\section{Preliminary}

For a vertex subset $S\subseteq V(G)$, the neighborhood of $S$ in $G$ is $N_G(S)=(\cup_{u\in S}N_G(u))\setminus S$; the induced subgraph $G[S]$ of $S$ in $G$ has vertex set $S$ and edge set $\{e=uv\in E(G): u,v\in S\}$. For  two vertex subsets $A,B\subseteq V(G)$, denote by $[A,B]_G$ the set of edges in $G$ with one end-vertex in $A$ and the other in $B$.

Let  $G_1$ and $G_2$ be two graphs. Define two natural projections $p_1$ and $p_2$ on $V(G_1)\times V(G_2)$ as follows: $p_1(x,y)=x$ and $p_2(x,y)=y$ for any $(x,y)\in V(G_1)\times V(G_2)$. Denote $\mathcal{G}=G_1\times G_2$. For any vertex $u\in V(G_1)$, let $G_2^{u}=\{(u,v)\in V(G_1\times G_2):v\in V(G_2)\}$ and call it the $G_2$-layer of $\mathcal{G}$ with respect to $u$. It is easy to see that  $[G_2^{u_{1}},G_2^{u_{2}}]_{\mathcal{G}}\neq\emptyset$  in $\mathcal{G}$  if and only if $u_{1}u_{2}\in E(G_1)$. Define $G_2^{S}=\cup_{u\in S}G_2^{u}$ for any vertex subset $S\subseteq V(G_1)$.

\begin{lem}
Let $A$ be a vertex subset of the graph $\mathcal{G}=K_2\times K_n$ ($n\geq5$). Assume $2\leq|A|\leq2n-2$. Then
$|[A,V(\mathcal{G})\setminus A]_{\mathcal{G}}|\geq 2(n-2)$, and the equality holds if and only if ($i$) $|A|=2$ and $\mathcal{G}[A]$ is isomorphic to $K_2$, or ($ii$) $|A|=2n-2$ and $\mathcal{G}[V(\mathcal{G})\setminus A]$ is isomorphic to $K_2$.
\end{lem}

\noindent{\bf Proof.} Let $U_1$ and $U_2$ be the bipartition of $\mathcal{G}$. Assume $A_i=A\cap U_i$ and $a_i=|A_i|$ for $i=1,2$. Without loss of generality, assume $|A|\leq |V(\mathcal{G})\setminus A|$, that is $|A|\leq n$. Then
\[
\begin{array}{ll}
|[A,V(\mathcal{G})\setminus A]_{\mathcal{G}}|&\geq a_1(n-1-a_2)+a_2(n-1-a_1)\\
&=(n-1)(a_1+a_2)-2a_1a_2\\
&\geq(n-1)|A|-2\frac{|A|^2}{4}\\
&\geq2(n-2).
\end{array}
\]
Furthermore, the equality holds if and only if $|A|=2$ and $\mathcal{G}[A]$ is isomorphic to $K_2$. $\square$

\begin{cor}
Assume $n\geq5$. Then ($i$) $K_{2}\times K_{n}$  is super-$\lambda$; ($ii$) $K_{2}\times K_{n}$  is super-$\lambda^{'}$.
\end{cor}

\noindent{\bf Proof.} Since $\lambda(K_{2}\times K_{n})\leq\delta(K_{2}\times K_{n})=n-1$ and $2(n-2)>n-1$, we obtain that $K_{2}\times K_{n}$  is super-$\lambda$ by Lemma 2.1.  For any vertex subset $A\subseteq V(K_{2}\times K_{n})$ with $3\leq|A|\leq 2n-3$, Lemma 2.1 shows that $|[A,V(K_{2}\times K_{n})\setminus A]_{K_{2}\times K_{n}}|>2(n-2)$. Thus, by $\lambda'(K_{2}\times K_{n})\leq\xi(K_{2}\times K_{n})=2(n-2)$, we prove that  $K_{2}\times K_{n}$  is super-$\lambda^{'}$. $\square$

By a similar argument as Lemma 2.1, we  have the following lemma.

\begin{lem}
Let $A$ be a vertex subset of the graph $\mathcal{G}=K_2\times T_n$ ($n\geq3$). Assume $2\leq|A|\leq2n-2$. Then
$|[A,V(\mathcal{G})\setminus A]_{\mathcal{G}}|\geq 2(n-1)$, and the equality holds if and only if ($i$) $|A|=2$ and $\mathcal{G}[A]$ is isomorphic to $K_2$, or ($ii$) $|A|=2n-2$ and $\mathcal{G}[V(\mathcal{G})\setminus A]$ is isomorphic to $K_2$.
\end{lem}

\begin{cor}
Assume $n\geq3$. Then ($i$) $K_{2}\times T_{n}$  is super-$\lambda$; ($ii$) $K_{2}\times T_{n}$  is super-$\lambda^{'}$.
\end{cor}

\section{Main Results}

\begin{thm}
For any nontrivial connected graph $G$ and any integer $n\geq5$. If $n(n-1)\lambda^{'}(G)>(n-1)\xi(G)+2(n-2)$, then $G\times K_{n}$ is super-$\lambda^{'}$.

\end{thm}

\noindent{\bf Proof.} Denote $\mathcal{G}=G\times K_{n}$. Since $n(n-1)\lambda^{'}(G)>(n-1)\xi(G)+2(n-2)$, we have $\lambda^{'}(\mathcal{G})=(n-1)\xi(G)+2(n-2)$ by Theorem 1.1.
Assume to the contrary that $\mathcal{G}$ is not super-$\lambda^{'}$.
Let $F$ be a minimum restricted edge-cut of $\mathcal{G}$. Then $|F|=(n-1)\xi(G)+2(n-2)$ and $\mathcal{G}-F$ has exactly two components $C_{1}$ and $C_{2}$, where
$|V(C_{1})|\geq3$ and $|V(C_{2})|\geq3$.

\noindent{\bf Case 1.} Each edge $e=(u_{1},v_{1})(u_{2},v_{2})\in E(C_{1})$ satisfies
$K_{n}^{u_{i}}\cap V(C_{2})=\emptyset$ for $i=1,2$.

Let $p_1(C_i)=U_i$ for $i=1,2$. By the assumption,  we have $V(C_i)=K_n^{U_i}$ for $i=1,2$. Then $[U_1, U_2]_G$ is a restricted edge-cut of $G$ and $|[U_1, U_2]_G|\geq\lambda'(G)$. Therefore
\[
|F|=|[V(C_{1}),V(C_{2})]_{\mathcal{G}}|=n(n-1)|[U_1, U_2]_G|\geq n(n-1)\lambda^{'}(G)>(n-1)\xi(G)+2(n-2),
\]
which is  a contradiction.

\noindent{\bf Case 2.} There is an edge $e=(u_{1},v_{1})(u_{2},v_{2})\in E(C_{1})$ such that $K_{n}^{u_{i}}\cap V(C_{2})\neq \emptyset$ for $i=1,2$.

Let $A_{1}=K_{n}^{u_{1}}\cap V(C_{1})$, $B_{1}=K_{n}^{u_{1}}\cap V(C_{2})$, $A_{2}=K_{n}^{u_{2}}\cap V(C_{1})$ and $B_{2}=K_{n}^{u_{2}}\cap V(C_{2})$.
Then $|A_{1}|+|B_{1}|=|A_{2}|+|B_{2}|=n$ and $|A_{i}|\geq1,|B_{i}|\geq1$ for $i=1,2$.
Since the induced subgraph $\mathcal{G}[K_n^{u_i}\cup K_n^{u}]$ is isomorphic to $K_2\times K_n$ for any neighbor $(u,v)$ of $A_i$ in $\mathcal{G}$ for $i=1,2$, we have $|[K_n^{u_i}, K_n^{u}]_{\mathcal{G}}\cap F|\geq n-1$ by Lemma 2.1.

Assume $|A_1|+|A_2|=2$. Since $|V(C_1)|\geq3$, there is a vertex $(u',v')\in V(C_1)$ such that $(u',v')$ is adjacent  to the vertex in $A_1$ or $A_2$. Without loss of generality, assume   $(u',v')$ is adjacent to the vertex in $A_1$. By Lemma 2.1, we obtain $|[K_n^{u_1}, K_n^{u_2}]_{\mathcal{G}}\cap F|\geq 2(n-2)$ and $|[K_n^{u_1}, K_n^{u'}]_{\mathcal{G}}\cap F|\geq 2(n-2)$. Thus
\[
\begin{array}{ll}
|F|&\geq\sum\limits_{u\in N_G(u_1)\setminus\{u_2,u'\}}|[K_n^{u_1}, K_n^{u}]_{\mathcal{G}}\cap F|
+\sum\limits_{u\in N_G(u_2)\setminus\{u_1\}}|[K_n^{u_2}, K_n^{u}]_{\mathcal{G}}\cap F|\\
&\ \ \ \ +|[K_n^{u_1}, K_n^{u_2}]_{\mathcal{G}}\cap F|+|[K_n^{u_1}, K_n^{u'}]_{\mathcal{G}}\cap F| \\
&\geq(n-1)(d(u_{1})-2)+(n-1)(d(u_{2})-1)+2(n-2)+2(n-2)\\
&=(n-1)(d(u_{1})+d(u_{2})-2)+2(n-2)+(n-3)\\
&\geq(n-1)\xi(G)+2(n-2)+(n-3)\\
&>(n-1)\xi(G)+2(n-2),
\end{array}
\]
which contradicts to the assumption. Similarly, the case $|B_1|+|B_2|=2$ can also obtain a contradiction.  So we assume $|A_1|+|A_2|\geq3$  and $|B_1|+|B_2|\geq3$. By Lemma 2.1, $|[K_n^{u_1}, K_n^{u_2}]_{\mathcal{G}}\cap F|>2(n-2)$. Thus
\[
\begin{array}{ll}
|F|&\geq\sum\limits_{u\in N_G(u_1)\setminus\{u_2\}}|[K_n^{u_1}, K_n^{u}]_{\mathcal{G}}\cap F|+\sum\limits_{u\in N_G(u_2)\setminus\{u_1\}}|[K_n^{u_2}, K_n^{u}]_{\mathcal{G}}\cap F|+|[K_n^{u_1}, K_n^{u_2}]_{\mathcal{G}}\cap F| \\
&>(n-1)(d(u_{1})-1)+(n-1)(d(u_{2})-1)+2(n-2)\\
&=(n-1)(d(u_{1})+d(u_{2})-2)+2(n-2)\\
&\geq(n-1)\xi(G)+2(n-2), \\
\end{array}
\]
which is a contradiction.

\noindent{\bf Case 3.} There is no edge $e'=(u',v')(u'',v'')\in E(C_{1})$ such that $K_{n}^{u'}\cap V(C_{2})\neq \emptyset$ and  $K_{n}^{u''}\cap V(C_{2})\neq \emptyset$. But there exists an edge $e=(u_{1},v_{1})(u_{2},v_{2})\in E(C_{1})$ such that $K_{n}^{u_1}\cap V(C_{2})\neq \emptyset$.

By this assumption, each vertex $w\in N_G(u_1)$ satisfies $K_{n}^{w}\subseteq V(C_{1})$ or
$K_{n}^{w}\subseteq V(C_{2})$. Let $M_1=\{x_1,\cdots,x_k\}$ be the neighbors of $u_1$ in $G$ such that $K_{n}^{x_i}\subseteq V(C_{1})$ for $i\in\{1,\cdots,k\}$; and let $M_2=\{y_1,\cdots,y_h\}$ be the neighbors of $u_1$ in $G$ such that $K_{n}^{y_j}\subseteq V(C_{2})$ for $j\in\{1,\cdots,h\}$,

\noindent{\bf Subcase 3.1.} There is a vertex $x\in N_G(M_1)$ such that $K_{n}^{x}\subseteq V(C_{1})$, and there is a vertex $y\in N_G(M_2)$ such that $K_{n}^{y}\subseteq V(C_{2})$.

Without loss of generality, assume $x\in N_G(x_1)$ and $y\in N_G(y_1)$. Then $K_{n}^{x}\cup K_{n}^{x_1}\subseteq V(C_{1})$ and $K_{n}^{y}\cup K_{n}^{y_1}\subseteq V(C_{2})$. Since there are at least $\lambda'(G)$ edge-disjoint paths between $\{x,x_1\}$ and $\{y,y_1\}$ in $G$, we obtain that there are at least $n(n-1)\lambda'(G)$ edge-disjoint paths between $K_{n}^{x}\cup K_{n}^{x_1}$ and $K_{n}^{y}\cup K_{n}^{y_1}$ in $\mathcal{G}$. Thus $|F|\geq n(n-1)\lambda'(G)>(n-1)\xi(G)+2(n-2)$, a contradiction.

\noindent{\bf Subcase 3.2.} Every vertex $x\in N_G(M_1)$ satisfies $K_{n}^{x}\cap V(C_{2})\neq\emptyset$, or every vertex $y\in N_G(M_2)$ satisfies $K_{n}^{y}\cap V(C_{1})\neq\emptyset$.

Assume, without loss of generality, that every vertex $x\in N_G(M_1)$ satisfies $K_{n}^{x}\cap V(C_{2})\neq\emptyset$.

If $|K_n^{u_1}\cap V(C_2)|\geq2$, then $|[K_n^{u_1}, K_n^{x_1}]_{\mathcal{G}}\cap F|>2(n-2)$ by Lemma 2.1. Thus
\[
\begin{array}{ll}
|F|&\geq\sum\limits_{u\in N_G(u_1)\setminus\{x_1\}}|[K_n^{u}, K_n^{u_1}]_{\mathcal{G}}\cap F|
+\sum\limits_{x\in N_G(x_1)\setminus\{u_1\}}|[K_n^{x}, K_n^{x_1}]_{\mathcal{G}}\cap F|+|[K_n^{u_1}, K_n^{x_1}]_{\mathcal{G}}\cap F| \\
&>(n-1)(d(u_1)-1)+(n-1)(d(x_{1})-1)+2(n-2)\\
&=(n-1)(d(u_1)+d(x_{1})-2)+2(n-2)\\
&\geq(n-1)\xi(G)+2(n-2), \\
\end{array}
\]
which is a contradiction.

If $|K_n^{u_1}\cap V(C_2)|=1$, then $|[K_n^{u_1}, K_n^{y_1}]_{\mathcal{G}}\cap F|>2(n-2)$ by Lemma 2.1. Thus
\[
\begin{array}{ll}
|F|&\geq\sum\limits_{u\in N_G(u_1)\setminus\{y_1\}}|[K_n^{u}, K_n^{u_1}]_{\mathcal{G}}\cap F|+\sum\limits_{x\in N_G(x_1)\setminus\{u_1\}}|[K_n^{x}, K_n^{x_1}]_{\mathcal{G}}\cap F|+|[K_n^{u_1}, K_n^{y_1}]_{\mathcal{G}}\cap F| \\
&>(n-1)(d(u_1)-1)+(n-1)(d(x_{1})-1)+2(n-2)\\
&=(n-1)(d(u_1)+d(x_{1})-2)+2(n-2)\\
&\geq(n-1)\xi(G)+2(n-2), \\
\end{array}
\]
which is also a contradiction.

The proof is thus complete.
$\square$

By Theorem 3.1, we obtain the following corollary immediately.

\begin{cor}
For any nontrivial graph $G$ and any integer $n\geq5$. If $G$ is $\lambda^{'}$-optimal, then $G\times K_{n}$ is super-$\lambda^{'}$.
\end{cor}

Similar to the proof of Theorem 3.1, we can obtain the following result by using Theorem 1.2 and Lemma 2.3.

\begin{thm}
For any nontrivial graph $G$ and any integer $n\geq3$. If $n^{2}\lambda^{'}(G)>n\xi(G)+2(n-1)$, then $G\times T_{n}$ is super-$\lambda^{'}$.
\end{thm}

\begin{cor}
For any nontrivial connected graph $G$ and any integer $n\geq3$. If $G$ is $\lambda^{'}$-optimal, then $G\times T_{n}$ is super-$\lambda^{'}$.
\end{cor}

\end{document}